\documentclass[12pt]{article}
\usepackage[dvips]{graphics}
\begin{document}

\title{Critical behavior in Ricci flow}

\author{David Garfinkle
\thanks {Email: garfinkl@oakland.edu}\\
Department of Physics, Oakland University\\
Rochester, Michigan 48309\\
\\
James Isenberg
\thanks {Email: jim@newton.uoregon.edu} \\
Department of Mathematics, University of Oregon\\
Eugene, OR}

\maketitle

\begin{abstract}
We use numerical techniques to study the formation of singularities 
in Ricci flow. 
Comparing the Ricci flows corresponding to a one parameter family 
of initial geometries on $S^3$ 
with varying amounts of $S^2$ neck pinching, we find critical 
behavior at the threshold of singularity
formation
\end{abstract}
\section{Introduction}

Given a manifold $M^n$ and a Riemannian metric ${\bar g}$ specified 
on $M^n$, the Ricci flow 
determines a one parameter family of metrics $g(t)$ via the 
geometric evolution equation
\begin{equation}
{\partial _t} {g_{ab}(t)}= - 2 {R_{ab}[g(t)]},
\label{ricflow}
\end{equation}
with initial condition 
\begin{equation}
g(0)={\bar g}.
\label{initdata}
\end{equation}
(Here $R_{ab}$ is the Ricci tensor of the metric $g_{ab}$.) Ricci flow 
has been a very effective tool for studying the sorts of special 
geometries which a manifold admits. In most of these applications to 
date \cite{Ham3} \cite{Ham2} \cite{Ham4}, the flow is shown to converge 
to the special geometries (e.g., to a constant negative curvature 
hyperbolic metric for the Ricci flow on any two dimensional manifold 
of genus greater than one). Future applications, however, are likely 
to require that one understand Ricci flows which develop singularities.

In this work, we begin a program of numerical study of Ricci flow 
singularities. We focus here on the following question: Say we have a 
one parameter family of metrics ${\bar g}_\lambda$ specified on the 
manifold $S^3$ and say we know that the Ricci flows $g_\lambda(t)$ 
starting at ${\bar g}_\lambda$ for very large $\lambda$  converge 
(with the volume suitably normalized) to the round sphere metric, 
while the Ricci flows $g_\lambda(t)$ for small values of the parameter 
$\lambda$ become singular, in the sense that (regardless of volume norm) 
the curvature of $g_\lambda (t)$ grows without bound as $t$ increases. 
What happens to the Ricci flows $g_\lambda (t)$ with intermediate 
values of $\lambda$? Is there, in particular, a certain threshold 
value $\lambda_{crit}$ for which the Ricci flow has interesting 
intermediate behavior: neither convergence, nor formation of a 
standard singularity?  

For solutions of Einstein's equations representing gravitational 
collapse, this sort of question has been studied 
extensively \cite{matt} \cite{critreview}. 
Very distinct threshold behavior has been found, with a remarkable degree 
of universality. That is, in examining a number of one parameter families 
of initial data for gravitational collapse, one finds qualitatively the 
same sort of discretely self-similar solution occurring for threshold 
initial data for all of these families. 

For Ricci flow, our work here is the first search for critical behavior. 
We have 
examined a particular one parameter family of spherically 
symmetric 
``corseted sphere'' geometries on $S^3$, with the parameter describing 
the degree of corseting at the equator, and therefore parametrizing the 
magnitude of the $S^2$ neck pinch curvature at the equator. We do indeed 
find numerically that for geometries with a small amount of corseting, 
the Ricci flow converges to the round sphere metric, while for geometries 
with a large amount of corseting, an $S^2$ neck pinch singularity occurs. 
Moreover, we find that there is a critical value of the parameter, dividing 
the two regimes.  Finally, our studies show that the Ricci flow for the 
geometry marked by the critical parameter value neither converges to the 
round sphere geometry nor forms an $S^2$ neck pinch singularity.  Instead, 
this flow approaches a ``javelin" geometry,  marked by curvature 
singularities at the poles, with roughly uniform curvature between the 
poles. This javelin geometry  corresponds to the ``type 3" singularity 
described by Hamilton \cite{HamForm} and discussed by Chow.\cite{Chow}

We describe in detail in Section 2 the corseted sphere geometries that 
we study here, and write out the detailed form of the Ricci flow equations 
for these geometries. Also in Section 2 we describe our numerical methods. 
We present our results in Section 3, noting the behavior of the Ricci 
flows for subcritical, supercritical,  and critical initial geometries. 
Concluding remarks appear in Section 4.

\section{Corseted Sphere Geometries \\ and Their Flow Equations}

The corseted sphere geometries and their flows are all represented by 
spherically symmetric metrics on $S^3$ of the form 
\begin{equation}
g = {e^{2X}}  \left ( {e^{- 2 W}}  d {\psi ^2}  +  {e^{2 W}}
 {\sin ^2} \psi  [ d {\theta ^2}  +  {\sin ^2} \theta  d {\phi ^2} ]
\right )  
\label{metric}
\end{equation}
Here $(\psi, \theta, \phi)$ are standard angular coordinates on 
the three sphere; spherical symmetry holds so long as we assume 
that the metric functions $X$ and $W$ are functions only of $\psi$. 

Smoothness of the metric at the poles, where $\psi$ takes the values $0$  
and $\pi$, requires that $X$ and $W$ be even functions of $\sin \psi$
in a neighborhood of the poles.  Therefore ${\partial _\psi} X$ and
${\partial _\psi} W$ must vanish at the poles. 
However, smoothness of the metric also requires that $W$ vanish at the 
poles.
To avoid the numerically awkward 
imposition of two conditions on $W$ at a single point, it is convenient 
to replace the variable $W$ by $S \equiv W/{\sin ^2} \psi$. Smoothness 
of the metric at the poles can then be enforced by the requirement that 
${\partial _\psi } X$ and ${\partial _\psi } S$ vanish 
at $\psi=0$ and $\psi=\pi$.

To obtain the corseted sphere geometries, we set $W=X$, and choose 
$X$ so that $ 4 {e^{4X}} {\sin ^2} \psi ={\sin ^2} 2 \psi $ 
for ${\cos ^2} \psi \ge 1/2$ and 
$4 {e^{4X}}{\sin ^2} \psi = {\sin ^2} 2 \psi + 4 \lambda {\cos ^2} 2 \psi $
for ${\cos ^2} \psi \le 1/2$. 
Here $\lambda$ is a constant, which parametrizes the degree of corseting for 
these geometries. For $\lambda=0$, the geometry represents two round three 
spheres joined at the poles. This is a singular geometry. For $\lambda$ 
positive, the cusp smooths out and the geometry is non singular; however, 
for small values of $\lambda$ one still expects the curvature at the tightly 
pinched equator to be very large. In Figure 1, we graph the area of the $\psi=constant$ cross-sections as a function of  $\psi$ for a few representative values of $\lambda$. 
\begin{figure}
\includegraphics{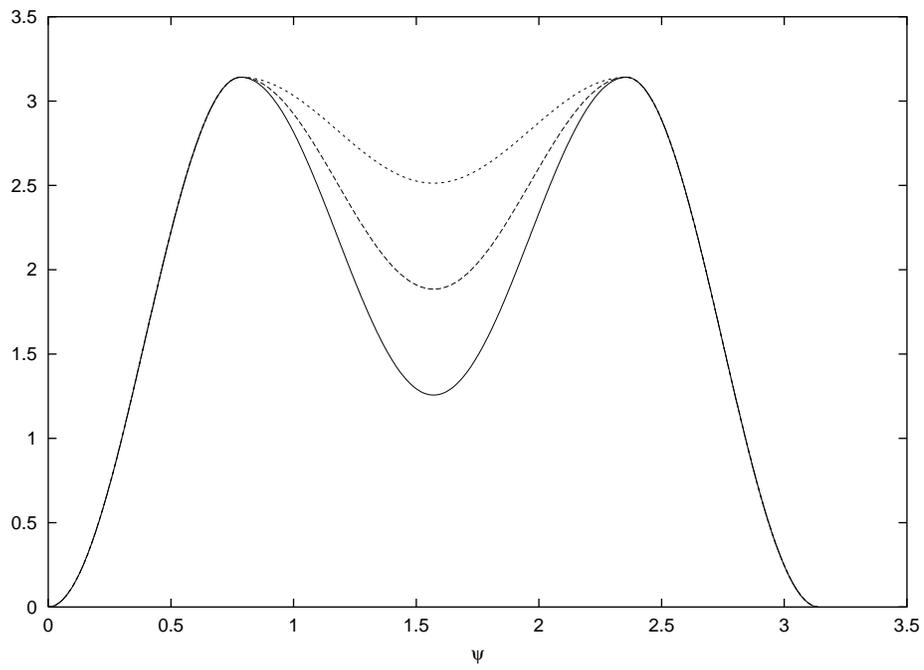}
\caption{\label{f1} area of the $S^2$s as a function of $\psi $ for 
$\lambda =0.1$ (solid line), $\lambda = 0.15$ (dashed line) and 
$\lambda = 0.2$ (dotted line).}
\end{figure}
In Figure 2, we graph the Ricci curvature eigenvalue in the 
direction along the $S^2$ symmetry as a function of $\psi$. Note that as 
a consequence of the assumption in these geometries that $W=X$, we 
verify that the coordinate $\psi$ gives the value of the geodesic 
distance in the radial direction. 
\begin{figure}
\includegraphics{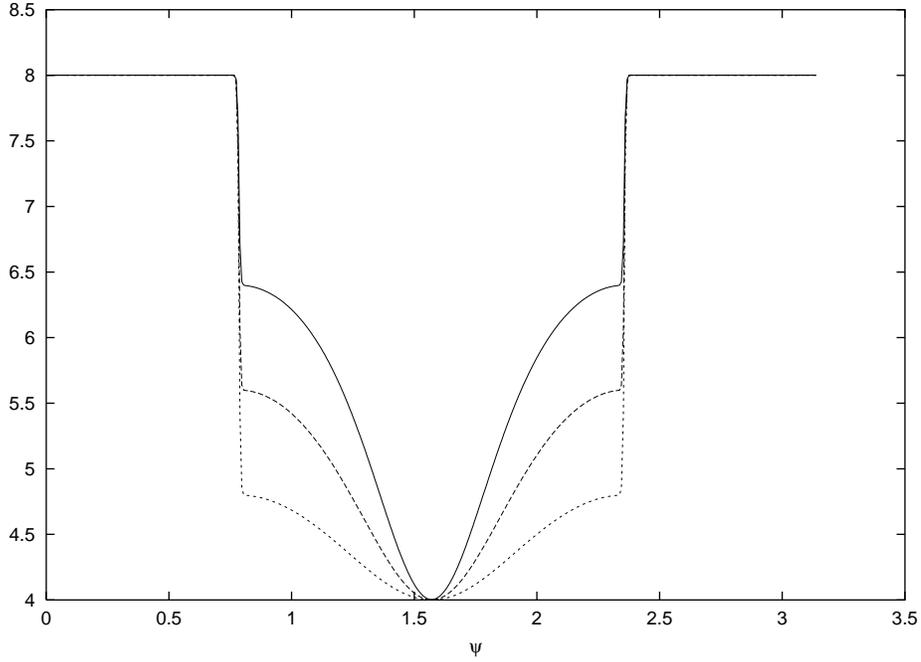}
\caption{\label{f2} $R_{S^2}$ as a function of $\psi $ for
$\lambda =0.1$ (solid line), $\lambda = 0.15$ (dashed line) and 
$\lambda = 0.2$ (dotted line).}
\end{figure}

We are interested in the Ricci flow of the corseted sphere metrics. 
However, since the Ricci flow equation (\ref{ricflow}) is only weakly 
parabolic, and since numerical evolutions appear to be more stable for 
strongly parabolic systems, we instead work with the DeTurck flow. 
\cite{DeT} The strongly parabolic PDE system generating the 
DeTurck flow is 
\begin{equation}
{\partial _t} {{\hat g}_{ab}} = -2{{\hat R}_{ab}}  + 2  {{\hat D}_{(a}}{V_{b)}} 
\label{DeTurck}
\end{equation}
where ${\hat D}_a$ is the derivative operator associated with
the metric ${\hat g}_{ab}$,  and where the vector field $V^a$ is given by
\begin{equation}
{V^a} = {{\hat g}^{bc}}  \left ( {{\hat \Gamma }^a _{bc}}  -  {\Delta ^a
_{bc}} \right )
\label{vector}
\end{equation}
with $ {{\hat \Gamma }^a _{bc}}$ being the connection of the metric
${\hat g}_{ab}$ and with $ {\Delta ^a _{bc}}$ being any fixed connection.
For a given initial geometry, one can find the corresponding Ricci flow 
$g(t)$ by first finding the corresponding DeTurck flow ${\hat g}_{ab}$, 
and then pulling back via the time dependent diffeomorphism generated by 
the vector field $V$ defined above. 

Note that neither the Ricci flow (\ref{ricflow})  nor the DeTurck flow 
(\ref{DeTurck}) preserves volume. One can normalize the volume for either 
of the flows by adding the term ${{2{\hat r}} \over 3} \; {{\hat g}_{ab}}$, 
with $\hat r$ the spatial average of the scalar curvature,  to the flow 
equation. Alternatively, one can control the volume along either flow by 
periodic uniform blowups.

We now calculate the volume normalized DeTurck flow for spherically 
symmetric metrics of the form (\ref{metric}). Choosing the reference 
connection to be that of the round sphere, and using primes to denote the 
spatial derivative $\partial_\psi$, we obtain the evolution equations
\begin{eqnarray}
{\partial _t} X = {e^{2(W-X)}}  \bigg [ {X''}  +  2  \cot \psi
{X' }  -  2  +  {1 \over 2}  ( {{[{X'}]} ^2}  +  {{[{W'}]} ^2} )  
+  3 
{X'}  {W'}
\nonumber
\\
+  (1 \, - \, {e^{- 4 W}})  \left ( {1 \over {2
{\sin ^2} \psi }}  +  1  +  2  \cot \psi  {W'} \right )  \bigg ]  +  
{{\hat r} \over 3}.
\label{evolveX}
\\
{\partial _t}W = {e^{2(W-X)}} \bigg [ {W''}  +  2  \cot \psi
{W'}  -  {1 \over 2}  ( {{[{X'}]}^2}  +  {{[{W'}]}^2} )  -  3 
{X'} \, {W'}
\nonumber
\\
+  (1  -  {e^{- 4 W}})  \left ( 1  -  {3 \over {2
{\sin ^2} \psi }}  -  2  \cot \psi  [{X'} + 2  {W'}] 
\right )  \bigg ].
\label{evolveW}
\end{eqnarray}
Here we note that the average scalar curvature $\hat r$ is given by 
\begin{equation}
{\hat r} = {2 \over N} \; {\int _0 ^\pi }  d \psi  {e^{X + 3 W}}  
\left ({e^{ - 4 W}}  - 1  -  4  \sin \psi  \cos \psi {W'}  +  
{\sin^2}
\psi [ 3 + {{({X'} + {W'})}^2} ] \right ),
\end{equation}
where the normalization constant $N$ is given by
\begin{equation}
N \equiv {\int _0 ^\pi } \; d \psi \; {e^{3X + W}} \; {\sin ^2} \psi.
\end{equation}
We also note that the only non vanishing component of the vector 
field V defined in (\ref{vector}) is
\begin{equation}
{V_\psi } = -1  \left ( 3  {W'}  +  {X'} + 2 \cot \psi  [ 1
- {e^{- 4 W}} ] \right )
\end{equation}

As discussed earlier, it is useful to work with the quantity  
$S \equiv W/{\sin ^2} \psi$, rather than $W$. Since this definition 
implies that  $ {W'}={\sin ^2} \psi {S'}+2 \sin \psi \cos \psi S$, we 
readily obtain
\begin{eqnarray}
{\partial _t} S = {e^{2(W-X)}} \biggl [ {S''} + 6 \cot \psi {S'}
- 8 S  - {3 \over {2 {\sin ^4} \psi }}
\left ( 1 - 4 W - {e^{-4W}}\right )  
\nonumber \\
+ {{1 - {e^{-4W}}}\over {{\sin ^2} \psi }}
\left ( 1 - 2 [\cot \psi {X'} + 2 \sin \psi \cos \psi {S'} + 4
{\cos ^2} \psi S]\right ) 
\nonumber \\
- {1 \over 2} \left ( {{[{X'}/\sin \psi ]}^2} + 
{{[\sin \psi {S'}+2 \cos \psi S]}^2} + 6 [{X'}/\sin \psi ] [ \sin \psi
{S'} + 2 \cos \psi S] \right ) \biggr ].     
\label{evolveS}
\end{eqnarray}
Equations (\ref{evolveX}) and (\ref{evolveS}) are the ones 
that we evolve numerically,
noting that wherever we encounter $W$ or ${\partial_\psi} W$ these are 
to be expressed in terms of $S$ and ${\partial_\psi} S$.

\section{Numerical Results}

To study these evolution equations numerically, we proceed as follows: 
We divide the spatial coordinate range $(0,\pi)$ of $\psi$ into $N-2$ 
pieces, so that we have  $\Delta \psi = \pi /(N-2)$. We choose $N$ grid 
points, including a pair which run outside the coordinate range. So the 
first spatial grid point is at $\psi = - \Delta\psi /2$, while the last 
is at $\psi = \pi + \Delta \psi /2$. Then a function of the form 
$F(\psi, t_0)$ for fixed time $t_0$ is replaced by a set of $N$ numbers 
${F_i}=F((i-{\textstyle {3 \over 2}})\Delta \psi , t_0)$ where 
$1\le i \le N$. Spatial derivatives are replaced by centered finite 
differences in the following way:
\begin{eqnarray}
{\partial_\psi F} \to {{{F_{i+1}}-{F_{i-1}}}\over {2 \Delta \psi}}
\\
{\partial_\psi \partial_\psi F} \to {{{F_{i+1}}+{F_{i-1}}-2{F_i}}\over
{{(\Delta \psi )}^2}}
\end{eqnarray}
For the time dependence of these functions, we choose a fixed time step $\Delta t$, and replace $F((i-{\textstyle {3 \over 2}})\Delta \psi , n\Delta t)$ by the numbers $F^n_i$.

Now, for an evolution equation of the form $\partial_t F(\psi,t)= G(\psi, t)$ we numerically evolve using 
the approximation
\begin{equation}
{F^{n+1}_i}={F^n_i}+\Delta t {G^n_i}
\end{equation}
This evolution is implemented for all values of $i$ except $1$ and $N$.
Note that these two ``ghost zones'' are not part of the manifold since
$\psi $ is not in the range $0 \le \psi \le \pi$.  At the 
ghost zones we use smoothness of the metric which implies that 
$X'$ and $S'$ vanish at the poles.  We implement
this condition  as ${X_1}={X_2}$ and ${X_N}={X_{N-1}}$ (and correspondingly
for $S$).  

Runs were done (on Unix and Linux workstations), starting from 
a wide range of initial corseted sphere geometries (parametrized by 
$\lambda$).  As noted above, for small $\lambda$ we 
expect that the geometry is sufficiently close to singular that the
distorted $S^3$ will pinch off into a singularity, while for sufficiently 
large $\lambda$ the flow may overcome the distortion to evolve the
data to a single round $S^3$.    

These expectations are confirmed by our numerical simulations.  Figures
3 and 4 show the results of a ``subcritical'' run, {\it i.e.} one 
which does not result in a singularity.  For this run we have chosen 
$\lambda=0.2  $.
Note that at late times we find that $X$ approaches a constant and
that $ S \to 0$.  These are the values for a round $S^3$.  

\begin{figure}
\includegraphics{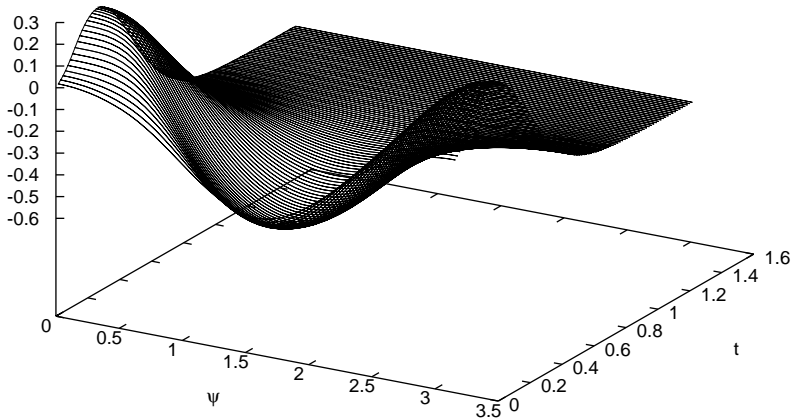}
\caption{\label{f3}$X$ for subcritical Ricci flow}
\end{figure}

\begin{figure}
\includegraphics{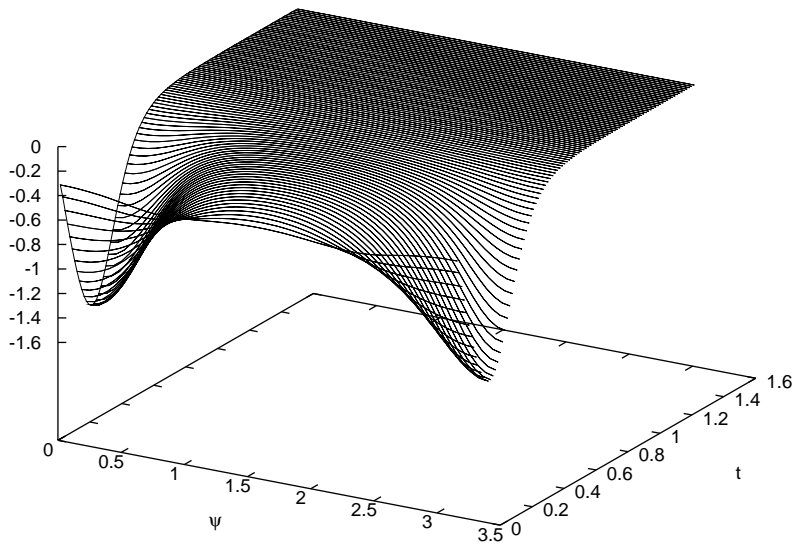}
\caption{\label{f4}$S$ for subcritical Ricci flow}
\end{figure}

Rather than focusing on the evolution of the metric components, 
it is more instructive to 
examine the behavior of curvature. As a consequence of the spherical 
symmetry of the geometries, the
Ricci tensor ${R^a}_b$ has two independent eigenvalues.  We call these
eigenvalues $R_{S^2}$ and $R_{\perp}$ where $R_{S^2}$ corresponds to 
the eigenspace
in the symmetry $S^2$ directions and $R_{\perp}$ corresponds to 
the eigenspace 
orthogonal to the symmetry $S^2$ directions.  Some straightforward
calculations show
\begin{eqnarray}
{R_{\perp}} = - 2 {e^{2(W-X)}} \left [ - 1 +  {X''} + {W''} +({X'} + 3 {W'})
\cot \psi + 2 ({X'}+{W'}){W'} \right ]    
\\
{R_{S^2}} = - {e^{2(W-X)}} \biggl [ - 2 + {{1 - {e^{-4W}}}\over {{\sin^2}
\psi}} + {X''} +{W''} + (3 {X'} + 5 {W'}) \cot \psi 
\nonumber \\
+({X'} + {W'})({X'} + 3 {W'})\biggr ]
\end{eqnarray}
(where we note that $W, {W'}$ and $W''$ can be expressed in terms of
$S, {S'}$ and $S''$).
One can express the invariants of the Ricci tensor (and since we are in 
3 dimensions, the invariants of the Riemann tensor as well) in terms of 
$R_{S^2}$ and $R_{\perp}$.  We have 
\begin{eqnarray}
R = 2{R_{S^2}} +  {R_{\perp}}
\\
{R^{ab}}{R_{ab}} =2 {R_{S^2}^2} + {R_{\perp} ^2} 
\\
{R^{abcd}}{R_{abcd}} = 4 {R^{ab}}{R_{ab}} - {R^2} = 2 {R_{\perp} ^2}
+ {{({R_{\perp}}-2{R_{S^2}})}^2}  
\end{eqnarray}
Figures 5 and 6 show the behavior of $R_{S^2}$ and $R_{\perp}$ 
for the same subcritical
run.  Note that both eigenvalues asymptotically approach the same 
constant value
at late times, which confirms the contention that the flow of 
a subcritical geometry converges to the round sphere geometry.

\begin{figure}
\includegraphics{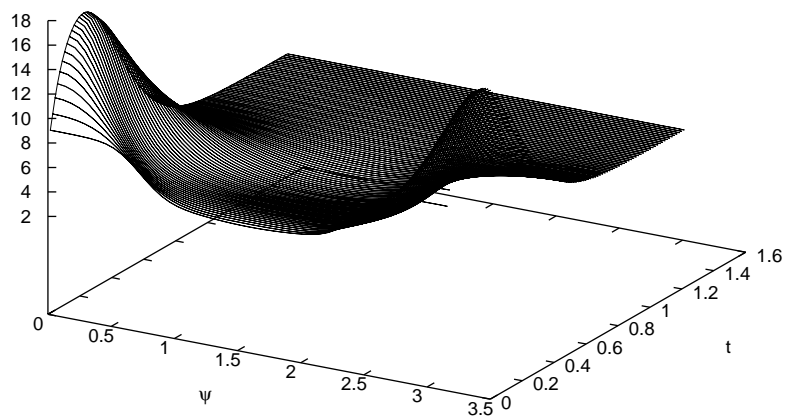}
\caption{\label{f5}$R_{S^2}$ for subcritical Ricci flow}
\end{figure}

\begin{figure}
\includegraphics{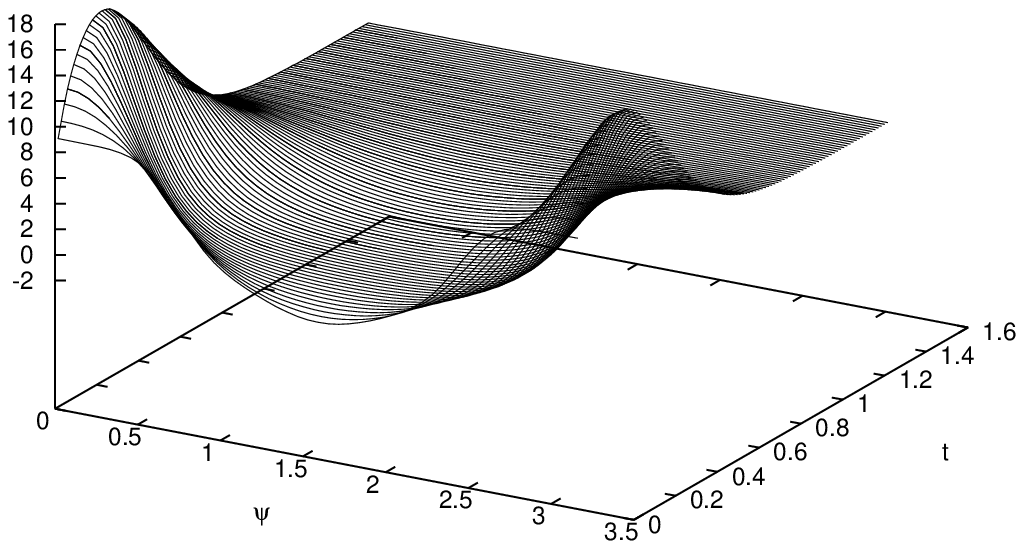}
\caption{\label{f6}$R_\perp $ for subcritical Ricci flow}
\end{figure}

We next consider the DeTurck flow of ``supercritical'' initial 
geometries, {\it i.e.} those for which the
evolution is singular.  Figures 7 and 8 plot the behavior of  $R_{S^2}$ 
and $R_{\perp}$ for an initial geometry with $\lambda=0.11 $. 
Characteristically, we find that in the neighborhood of the equator, 
$R_{S^2}$ grows without bound as $t$ increases, while $R_{\perp}$ 
appears to stay
bounded. As well, both are bounded away from the equator. 
This signals the formation of an $S^2$ neck pinch singularity at the equator.

\begin{figure}
\includegraphics{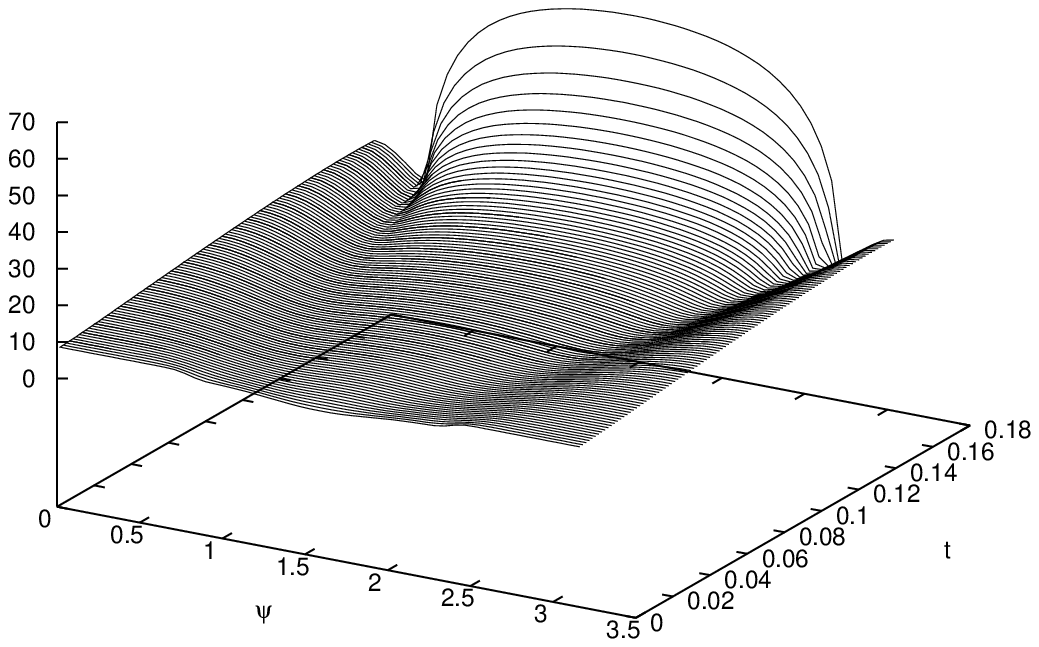}
\caption{\label{f7}$R_{S^2}$ for supercritical Ricci flow}
\end{figure}

\begin{figure}
\includegraphics{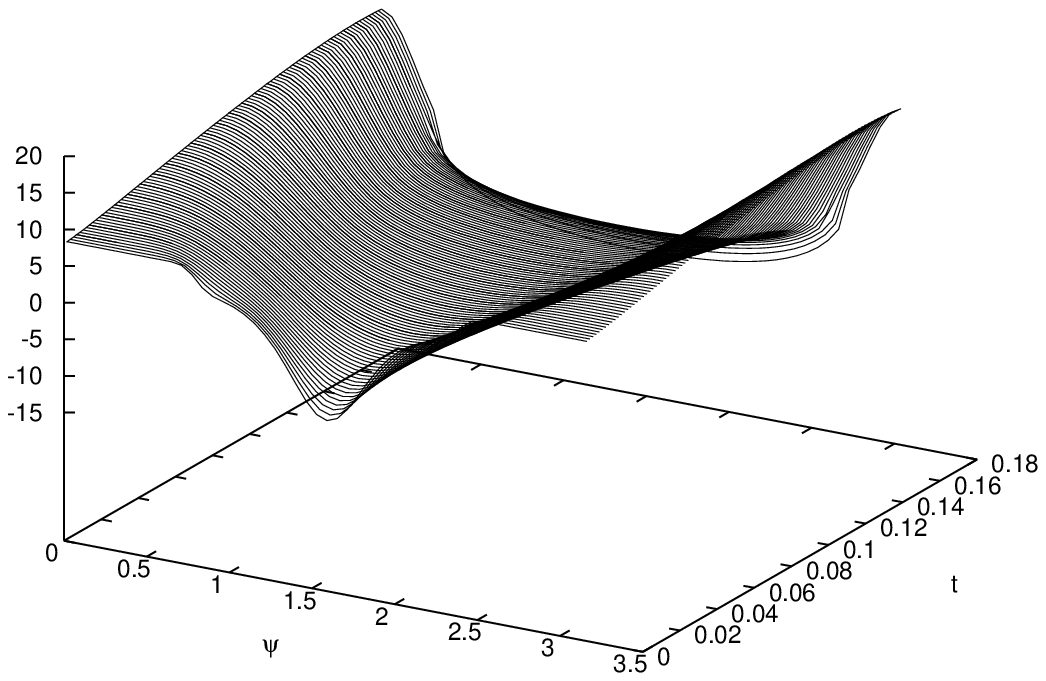}
\caption{\label{f8}$R_\perp$ for supercritical Ricci flow}
\end{figure}

The behavior just described for the DeTurck flow--and consequently for 
the 
Ricci flow--for $\lambda=0.11 $ corseted sphere data is also found for 
any initial geometry with $\lambda< 0.11$. Similarly, for $\lambda> 0.2$, 
the flow has the subcritical behavior illustrated in Figures 5 and 6. 
This has led us to seek critical behavior at a threshold value. Using a 
binary search, we have located the threshold value at approximately 
$\lambda=0.1639$

\begin{figure}
\includegraphics{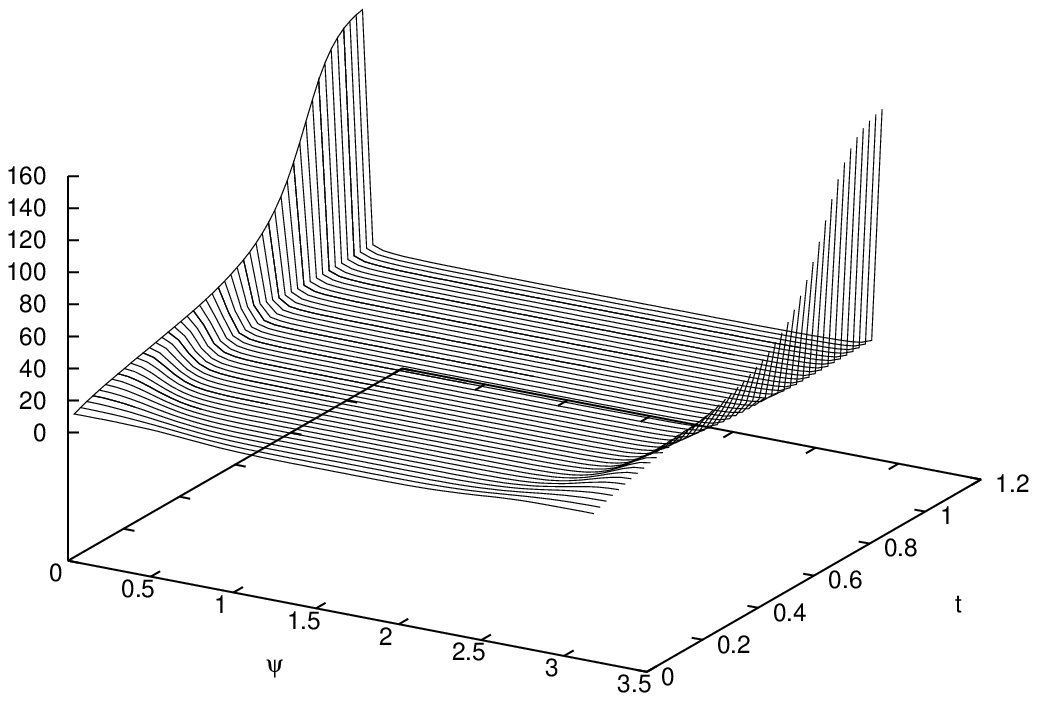}
\caption{\label{f9}$R_{S^2}$ for critical Ricci flow}
\end{figure}

\begin{figure}
\includegraphics{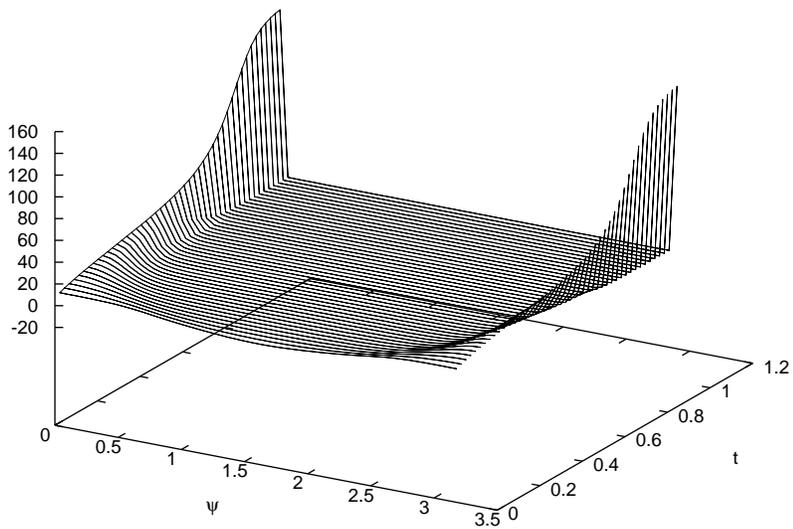}
\caption{\label{f10}$R_\perp$ for critical Ricci flow}
\end{figure}

The DeTurck flow starting at a corseted sphere geometry with this value of 
$\lambda$ behaves differently than the flow for both subcritical and 
supercritical geometries. As seen in Figures 9 and 10, $R_{\perp}$ gets 
small everywhere except at the poles, while $R_{S^2}$ slowly grows every 
where except at the poles. At the poles, both curvatures get very large. 
In a sense, the geometry approaches that of a three dimensional javelin, 
with $S^2$ cross-sections.  

\section{Conclusions}

The numerical work we have described here clearly shows rather special 
behavior for the Ricci flow at the threshold parameter value for a one 
parameter family of corseted spheres. Is this behavior in any sense 
universal? This has not yet been determined. We plan to examine other 
families of initial geometries to see if the threshold behavior persists. 
One interesting set of geometries we plan to examine are those with, 
initially, more than one neck pinch.. Do the neck pinches coalesce? 
Do we get javelin geometries for threshold initial data? We also plan to 
consider families of initial data which are not spherically symmetric.

\section{Acknowledgments}

This work was partially supported by NSF grants
PHY-9988790 to Oakland University and PHY-0099373 to The University of 
Oregon. We also thank the University of California at San Diego for 
hospitality while some of this work was carried out.

\end{document}